\documentclass[12pt]{amsart}
\usepackage{xypic} 
\LaTeXdiagrams          

\usepackage{amsthm,amsfonts, amssymb, amscd}
\newtheorem{prop}{Proposition}[section]
\newtheorem{def-prop}{Definition-Proposition}[section]
\newtheorem{thm}{Theorem}[section]
\newtheorem{cor}{Corollary}[section]

\newtheorem{lemma}{Lemma}[section]

\theoremstyle{remark}
\newtheorem{remark}{Remark}[section]
\newtheorem{example}{Example}[section]
\theoremstyle{definition}
\newtheorem{defn}{Definition}[section]
\def\P{{\bf P}}

\def\Z{{\mathbb Z}}
\def\bG{{\mathbb G}}
\def\endproof{\qed}

\DeclareMathOperator{\adj}{Adj}
\DeclareMathOperator{\im}{im}
\DeclareMathOperator{\td}{Td}

\def\Q{{\mathbb Q}}

\def\V{{\mathcal V}}
\def\O{{\mathcal O}}

\def\F{{\mathcal F}}
\def\t{\tau}

\def\Ibar{\overline{I}}
\def\a{{\mathfrak a}}

\def\liminv{{\displaystyle \lim_{\leftarrow (V,U) \in {\mathcal V}}}}
\begin{document}

\title[Equivariant Riemann-Roch]
{Riemann-Roch for equivariant Chow groups}
\author{Dan Edidin}
\address{Department of Mathematics\\ University of Missouri\\
Columbia MO 65211}
\email{edidin@math.missouri.edu}
\author{William Graham}
\address{Department of Mathematics\\ University of Georgia\\
Boyd Graduate Studies Research Center\\Athens, GA 30602}
\email{wag@math.uga.edu}
\thanks{The first author was partially by NSA grant
MDA904-97-1-0030
and the M.U. research board. The second author
was partially supported by the NSF}
\maketitle
\section{Introduction}
The purpose of this paper is to prove an equivariant
Riemann-Roch theorem for schemes or algebraic
spaces with an action
of a linear algebraic group $G$.  For a $G$-space $X$, this
theorem gives an
isomorphism
$$
\tau^G: G^G(X) \rightarrow \widehat{G^G(X)}_{\Q} \stackrel{\simeq}
\rightarrow \prod_{i=0}^{\infty} CH^i_G(X)_{\Q}.
$$
Here 
$\widehat{G^G(X)}$ is 
the completion of the equivariant Grothendieck group of coherent
sheaves along the augmentation ideal
of the representation ring $R(G)$, and
the groups $CH^i_G(X)$ are the equivariant Chow groups
defined in \cite{EIT}.  The map $\tau^G$ has the same
functorial properties as the non-equivariant Riemann-Roch map
of \cite{BFM}, \cite[Theorem 18.3]{Fulton} and 
if $G$ acts freely
then $\tau^G$ can be identified with the non-equivariant Todd
class map
$\tau_{X/G}: G(X/G) \rightarrow  CH^*(X/G)_{\Q}$.

The key to proving this isomorphism is a geometric description of
completions of the equivariant Grothendieck group (Theorem \ref{t.comp}).
Besides Riemann-Roch, this result has 
some purely $K$-theoretic applications.  In
particular, we prove (Corollary \ref{c.kock})
a conjecture of K\"ock 
(in the case of regular
schemes over fields) and extend to arbitrary characteristic a result of Segal on
representation rings (Corollary \ref{c.segal}).

For actions with finite stabilizers the equivariant Riemann-Roch
theorem is more precise; it gives
an isomorphism between a localization of $G^G(X)_\Q$ and
$\oplus CH^i_G(X)_\Q$ (Corollary \ref{c.fs}).  This formulation enables us
to give a simple proof of a conjecture of Vistoli (Corollary \ref{c.vistoli}).
If $G$ is diagonalizable, then we can express
$G^G(X)$ in terms of the equivariant Chow groups (an unpublished
result of Vistoli, cf. also \cite{Toen}).  
Actions with finite stabilizers are particularly important because
quotients by these actions arise naturally in geometric invariant
theory. 
In a subsequent paper,
we will use these results to express the Todd class map
for a quotient of such an action in terms of
equivariant Todd class maps, generalizing Riemann-Roch formulas of
Atiyah and Kawasaki.

The main tool of this paper is the approximation of total
space of the classifying bundle $EG$ by an open
subset $U$ of a representation $V$, where $G$ acts freely on $U$, and
$V-U$ is a finite union of {\it linear subspaces}.  Approximations to
$EG$ by open sets in representations were introduced by Totaro in Chow
theory \cite{Totaro}, and used in \cite{EIT} to define equivariant
Chow groups. However, in these papers, $V-U$ is 
only required to have large codimension: because
Chow groups are naturally graded we can identify
$\oplus_{i=0}^{N}CH^i_G(X)$ with $\oplus_{i=0}^{N}CH^i_G(X \times U)$ as
long as $\mbox{codim }(V-U)>N$. Since Grothendieck groups are not
naturally graded we
need the stronger condition that $V-U$ is a union of linear subspaces
to compare $G^G(X)$ with $G^G(X \times U)$.

Such $V$ and $U$ can be found for
tori or for Borel subgroups of $GL_n$.  For these groups, we prove that
the completion of $G^G(X)$ along the augmentation ideal has a geometric
description. This fact directly implies the equivariant Riemann-Roch
isomorphism for these groups.  For general $G$, however, it seems
unlikely that such $V$ and $U$ exist, so we must employ a less direct
approach.  The Riemann-Roch theorem for general $G$ is deduced by
embedding $G$ into $GL_n$ and then reducing the case of $GL_n$ to that
of a Borel subgroup. This strategy of proof is due to Atiyah and Segal for
compact groups. 

The necessity of using completions of equivariant Grothendieck groups
also goes back to Atiyah and Segal.  In our setting it is motivated as
follows.  For smooth varieties the Todd class map is defined by a
power series which, in the non-equivariant case, terminates on
any particular variety.  However, the equivariant Chow groups of a
fixed variety can be nonzero in arbitrarily high degree, so the
natural target of the equivariant Todd class map is 
the infinite product $\prod_{i=0}^{\infty} CH^i_G(X)_\Q$.
To obtain a Riemann-Roch isomorphism, it is natural to expect that the
equivariant Grothendieck groups must also be completed, as is indeed
the case.  There is one essential difference between these
completions.  The completion map for Chow groups is injective,
as it simply replaces a direct sum
by a direct product, but on the Grothendieck group side, information
is definitely lost by completing (cf. Section 5).

We thank the referee who read the paper very carefully and
suggested many corrections and improvements.

\subsection{Contents}
In Section 2 we define the
completions we need for the main results, and give a geometric
description of these completions for tori and for the subgroup of
upper triangular matrices of $GL_n$.  In Section 3 we construct the
equivariant Riemann-Roch map, and prove that it has the same functorial
properties as the non-equivariant Riemann-Roch map.  
Moreover, it behaves naturally with respect to
restriction to a subgroup (Section \ref{s.changegroups}). In Section \ref{s.weyl}
we illustrate the
use of this theorem by deriving the Weyl character formula for $SL_2$
following \cite{Bott}.  In Section 4 we prove that the Riemann-Roch map induces an
isomorphism on the completions.  Section 5 contains results, mentioned
above, on actions with finite stabilizers, and proves Vistoli's
conjecture.  Finally, in Section 6 we prove that two naturally
defined completions of equivariant $K$-theory are equal, 
and apply it to prove K\"ock's
conjecture and the result about representation rings.

\subsection{Conventions and notation}
All groups in this paper are assumed to be linear algebraic
groups over an arbitrary field $k$ (i.e. closed subgroup-schemes of
$GL_n(k)$).  Note that in characteristic $p$ such groups
need not be smooth over $k$. 
Representations of $G$ are assumed to be rational,
i.e., linear actions of $G$ on finite dimensional $k$-vector spaces. 

\medskip

{\bf Free actions.}
By a free action of $G$ on $X$ we mean an action that is
scheme-theoretically free, i.e., the action map
$G \times X \rightarrow X \times X$ is a closed embedding.
An action that is proper and set-theoretically free is free 
\cite[Lemma 8]{EIT}.

\medskip

{\bf Algebraic spaces.} This paper is written in the language of
algebraic spaces.  Unless stated otherwise, a space is an
equidimensional quasi-separated algebraic space of finite type over
$k$.  One reason to work in this category is that we need quotients of
the form $X \times^G U$ where $U$ is an open set in a representation
of $G$ on which $G$ acts freely.  In the category of algebraic spaces,
such quotients exist, by a result of Artin \cite[Proposition
22]{EIT}.  It is possible to work entirely in the smaller category of
schemes of finite type over $k$, but some mild technical hypotheses
are required to ensure that quotients exist as schemes
\cite[Proposition 23]{EIT}.  Because free actions are proper, if $X$
is separated, so is $X \times^G U$ \cite[Corollary 2.2]{EdMo}.

\medskip

{\bf Equivariant Chow groups}
We will use the notation $A_k^G(X)$ to refer to the equivariant Chow
groups of ``dimension'' $k$ as defined in \cite{EIT}. However, it is
usually more convenient to index by codimension.  We write $CH^i_G(X)$
for the ``codimension'' $i$ equivariant Chow group: if $X$ has pure
dimension $n$, then $CH^i_G(X) := A_{n-i}^G(X)$.  The notation
$A^i_G(X)$ refers to the degree $i$ operational Chow group;
when $X$ is smooth, $CH^i_G(X) \simeq A^i_G(X)$ \cite{EIT}. The
operational Chow groups have a natural product, and
$\oplus_{i=0}^{\infty}CH^i_G(X)$ is a positively graded module for the
graded ring $\oplus A^i_G(X)$.

The $G$-equivariant Chow ring of a point will be denoted
$A^*_G$. Pullback from a point makes $\oplus_i^\infty CH^i_G(X)$
a graded $A^*_G$-module.

\medskip

{\bf Equivariant $K$-theory}
We will use the notation $K^G(X)$ to refer to the
Grothendieck group of $G$-equivariant vector bundles.
The Grothendieck group of $G$-equivariant coherent sheaves
will be denoted $G^G(X)$. Tensor product gives
a ring structure on $K^G(X)$, and $G^G(X)$ is 
a module for this ring. If $X$ is regular and has
the property that every coherent sheaf is the quotient
of a locally free sheaf (this holds for regular
schemes, but is not known in general for regular algebraic spaces)
then a result of Thomason
\cite[Corollary 5.8]{Tho87} says that these groups are
isomorphic.

The representation ring $R(G)$ can be identified with $K^G(pt)$,
and consequently $G^G(X)$ has an $R(G)$-module structure.
If $A$ is a ring define $G^G(X)_A = G^G(X) \otimes_{\Z} A$, the
Grothendieck group with coefficients in $A$, and similarly
for Chow groups.

\section{Completions} \label{s.completions}
In this section (except in Example \ref{e.noncommute}) we work with a
fixed coefficient ring $A$ and write simply $G^G(X)$ for $G^G(X)_A$,
and similarly for Chow groups.

We will use the following completions in the sequel.
The graded ring $\oplus A^i_G(X)$ completes to
$\prod_{i=0}^{\infty}A^i_G(X)$ and
$\prod_{i=0}^{\infty} CH^i_G(X)$
is naturally a module for this ring.  We will let $J$ denote
the augmentation ideal in $A^*_G$, and let
$\widehat{CH^*_G(X)}$ denote the completion of the
$A^*_G$-module $CH^*_G(X)$ along $J$.

Completing equivariant $K$-theory requires more care since
$K^G$ is not naturally graded.  We complete as follows.

Choose an embedding $G \hookrightarrow GL_n$. This induces
a homomorphism $R=R(GL_n) \rightarrow R(G)$.  Thus
$G^G(X)$ is an $R$-module.  Let
$I$ be the augmentation ideal of $R$ (the ideal of
virtual representations of $GL_n$ of dimension $0$), and $I_G$ the
augmentation ideal of $R(G)$.  We denote
the $I$-adic completion of $G^G(X)$ by $\widehat{G^G(X)}$.    

This definition is convenient for proving results about $G$ by
reducing to the case of $GL_n$, but it is not obvious that it is
independent of the choice of embedding of $G$ into $GL_n$.  However,
this follows from Corollary \ref{c.segal}, which implies that $I$ and
$I_G$ give the same topology on $R(G)$ whenever $G \hookrightarrow
GL_n$; hence $I$-adic and $I_G$-adic completions are isomorphic.  A
special case of this result is proved in Proposition \ref{p.borelsegal}.  
In characteristic 0,  Corollary \ref{c.segal} is the
same as \cite[Corollary 3.9]{Segal}.  However, his methods do not
extend to characteristic $p$.

\subsection{Geometric completions}
There is a way of completing equivariant $K$-groups and Chow
groups that is more directly related to the definition of
equivariant Chow groups.  
For any representation
$V$ of $G$, let $V^0$ denote the open subset of points of $V$
on which $G$ acts freely and let $U$ be an open subset of $V^0$.  A
quotient $(X \times U)/G$ exists; such a quotient is usually written
as $X \times ^G U$.  Choose $V$ and $U$ such that the codimension of $V -
U$ is greater than $k$ (this is always possible by \cite{EGchar}); then by
definition, $CH^k_G(X) = CH^k(X \times^G U)$.

With this as motivation, let ${\mathcal V}$ be a collection of 
pairs  $(V,U)$ of $G$-modules and invariant open
sets with the following properties:
\medskip

(i) $G$ acts freely on $U$.\\

(ii) If $(V_1,U_1)$ and $(V_2,U_2)$ are in ${\mathcal V}$
then there is a pair $(V_1 \oplus V_2, U)$ in ${\mathcal V}$
such that $U$ contains $U_1 \times V_2$ and $V_1 \times U_2$.
\medskip

The elements of such a system can be partially
ordered by the rule $(W,U_W) < (V,U_V)$ if the $G$-module $V$ 
can be written as a direct sum
$V = W \oplus W'$ with $U_V \supset U_W \times W'$.
Suppose that $(V,U)$ is in ${\mathcal V}$;
by the homotopy property, we can identify $G^G(X) = G^G(X \times V)$ \cite
[Theorem 4.1]{Tho87} and 
$CH^i_G(X \times V)= CH^i_G(X)$ \cite{EIT}. In this way we obtain
surjective maps
$$k_V:G^G(X) \rightarrow G^G(X \times U)$$
and
$$r_V:\oplus CH^i_G(X) \rightarrow \oplus CH^i_G(X \times U).$$

\begin{lemma} \label{l.contain}
If $(W,U_W) < (V,U_W)$ then 
$\ker k_W \supseteq \ker k_V$ and
$\ker r_W \supseteq \ker r_V$.
\end{lemma}

\begin{proof}
Write $V = W \oplus W'$.  Then $U_V \supset U_W \times W'$,
so we have surjective maps
$$
G^G(X) \rightarrow G^G(X \times U_V) \rightarrow
G^G(X \times U_W \times W') \cong G^G(X \times U_W).
$$
The first map is $k_V$ and the composition
is $k_W$, proving the first inclusion.
A similar argument works for $r_V$ and $r_W$.
\end{proof}

\medskip

For $(W,U_W) < (V,U_V)$, the lemma implies there is a natural surjective
map
$$
G^G(X)/\ker k_V  \rightarrow G^G(X)/\ker k_W
$$
making these groups into an inverse system indexed by 
${\mathcal V}$.  A similar statement holds for the
Chow groups. Taking the inverse limit we obtain a completion
of $G^G(X)$. For an arbitrary system ${\mathcal V}$
it is difficult to describe this completion. However,
there is one situation where we can understand it.

Call a system ${\mathcal V}$ with properties (i) and (ii)
above {\it good} if it has a third property:

\medskip

(iii) If $(V,U)$ is in ${\mathcal V}$ then $V-U$ is contained
in a finite union of invariant linear subspaces.

\begin{thm} \label{t.comp}
Let ${\mathcal V}$ be a good system of representations.
Then:

(a) The topology on $G^G(X)$ induced by the subgroups
$\ker k_V$ coincides with the $I_G$-adic topology.  Hence
$\liminv (G^G(X)/\ker k_V)$ is isomorphic to the $I_G$-adic
completion of $G^G(X)$.

(b) The topology on $CH^*_G(X)$ induced by the subgroups
$\ker r_V$ coincides with the $J$-adic topology.  Hence
$\liminv CH^*_G(X)/\ker r_V$ is isomorphic to the
$J$-adic completion of $CH^*_G(X)$.

\end{thm}

The following result shows that for certain classes of
groups, good systems exist.

\begin{thm} \label{t.good}
If $G$ is a torus
or more generally any subgroup of the group of upper triangular
matrices in $GL_n$, then there exists a system ${\mathcal V}$ of good 
pairs for $G$.

In particular, if the ground field is algebraically closed
then any connected solvable group has a good system.
\end{thm}

The proofs of these theorems will be given in the next 2 subsections.
We close with two propositions which will be used
in the sequel.
\begin{prop} \label{p.brion}
Let $G$ be a connected subgroup of the group of upper triangular matrices.  
Then 
$$\widehat{CH_G^*(X)} \simeq \prod_{i=0}^{\infty} CH^i_G(X).$$
\end{prop}
\begin{proof}[Proof of Proposition \ref{p.brion}]
If $G$ is a connected subgroup of the upper triangular matrices
then $G = TU$ with $T$ a torus and $U$ unipotent 
\cite[Theorem 10.6]{Borel}.
An argument similar to that used in Thomason \cite[Proof of Theorem 1.11]
{Thoduke88}
implies that $CH^*_G(X) = CH^*_T(X)$ so we assume
$G = T$ is a torus.

Let $M = CH^*_T(X)$ and
$M_n = \oplus_{i = n+1}^{\infty} CH^i_T(X)$.  We will show
that the filtrations of $M$ given by $\{ J^n M \}$ and $\{ M_n \}$
give the same topology.

{\it Step 1.} For any $n$ we must show there exists $l$ such that $J^lM
\subseteq M_n$.  Since multiplication by $J^l$ maps $CH^i_T(X)$ to $M_{i+l-1}$,
if $l \geq n+1$ the desired inclusion holds.

{\it Step 2.} For any $n$ we must show there exists $l$ such that
$M_l \subseteq J^nM$.  Brion \cite[Theorem 2.1]{Brion} has shown that 
$CH^*_T(X)$ is generated as an $A^*_T$-module by fundamental
classes of $T$-invariant subvarieties.  Any such fundamental class
lies in $CH^i_T(X)$ for some $i \leq N$, where $N =\mbox{dim }X$.  
Hence if
$r > \mbox{dim }X$, then $CH^r_T(X) \subseteq \oplus_{i \leq N} A^{r -i}_T
\cdot CH^i_G(X) $.  Since $T$ is a torus
$A^*_T$ is generated as in algebra in degree 1 \cite[Section 3.2]{EIT}.
Thus, $A_T^{r-i} \subset J^{r-i}$ and we can conclude that if $l > n + N$,
we have $M_l \subseteq J^nM$.
\end{proof}

\begin{prop} \label{p.borelsegal} If $G =GL_n$ and $B$ (resp. $T$) is the group
of upper triangular (resp. diagonal) matrices then\\ (a) The map $R(B)
\rightarrow R(G) = R(B)^W$ is finite.\\ (b) $I_{GL_n}$ generates the same
topology on $R(B)$ (resp. $R(T)$) as $I_B$ (resp. $I_T$).
\end{prop}

\begin{proof} The Weyl group of $GL_n$ is $S_n$ over any
field, so by \cite[Example 3.8]{Serre}
$R(G) = R(B)^{S_n}$. 
In this case
it is well known that $R(B) \simeq \Z[t_1, \ldots , t_n, t_1^{-1},
\ldots , t_n^{-1}]$ and $R(G) = R(B)^{S_n} \simeq \Z[e_1, \ldots , e_n,
e_1^{-1}, \ldots e^{-1}_n]$ where $e_i$ is the $i$-th elementary
symmetric polynomial in the variables $t_1, \ldots , t_n$. This
fact obviously implies (a). Part (b) follows from
\cite[Corollary 3.9]{Segal}
applied to the maximal compact subgroups $U_n \subset GL_n({\bf C})$ and
$(S^1)^n \subset B$.
\end{proof}

\subsection{Proof of Theorem \ref{t.good}}


Observe that if $H \subset G$ is a closed subgroup and if
$\V$ is a good system for $G$, then it is also a good
system for $H$.
Thus, it suffices
to construct a good system for the group $B$ of upper
triangular matrices in $GL_n$.
We will do this 
as follows: Let $V_1$ be the vector space of upper triangular
matrices; $B$ acts on $V_1$ by left matrix multiplication.  Let $U_1
\subset V_1$ be the subset of invertible elements. Then 
$V_1 - U_1$ is a union of $n$ hyperplanes.  Set
$V_k = V_1^{\oplus k}$ and let $U_k \subset V_k$ consist of the
$k$-tuples of $V_1$ such that at least one element of the $k$-tuple
lies in $U_1$.  Then $V_k - U_k$ is a union of linear subspaces. By construction
the collection of pairs $\{(V_k,U_k)\}$ satisfies properties (ii) and
(iii) 
of a good system. The action of $B$ on $U_k$ has trivial
stabilizers, but to verify (i) we must check
that the action map $B \times U_k \rightarrow U_k \times U_k$
is a closed embedding.
A tuple of matrices $(A_1, \ldots , A_k, B_1, \ldots , B_k) \in
U_k \times U_k$ is in the image of $B \times U_k$ if and
only if there is a matrix $A \in B$ such that
$B_i = A A_i$ for all $i$. Thus, 
the image of
$B \times U_k$ in $U_k \times U_k$ is the closed
subvariety defined by the equations
$$B_i\adj(A_i) A_j = \det(A_i) B_j.$$

\begin{remark}
Let $(V,U)$ be a pair in our system $\V$ of good representations.
In the proof of Theorem \ref{t.comp} below, we only use 
the fact that $G$ acts with trivial stabilizers
(as opposed to freely) on the open set $U \subset V$.
However, the freeness of the action will be essential
when we apply Theorem \ref{t.comp} to prove the Riemann-Roch isomorphism.
The reason is that we will need to know that if $X$ is a
separated algebraic space then $X \times^G U$ is still
separated.
\end{remark}

\subsection{Proof of Theorem \ref{t.comp}}

\begin{lemma}\label{component} Let $X$ be a $G$-space which is
a union of invariant irreducible components $X_1,\ldots, X_k$.
Then proper push-forward gives a surjection 
$$ \oplus G^G(X_i) \rightarrow G^G(X). $$

\end{lemma}

\begin{proof} 
By induction it suffices to prove the lemma when $X = X_1 \cup X_2$.
Let $\tilde{X}$ be the disjoint union of $X_1$ and $X_2$.
Then $G^G(\tilde{X}) = G^G(X_1) \oplus G^G(X_2)$.
The finite surjective map $\tilde{X} \rightarrow X$
gives a map of localization exact sequences
$$\def\objectstyle{\scriptstyle}\def\labelstyle{\scriptstyle}
\xymatrix{
G^G(X_1 \cap X_2)^{\oplus 2} \ar[d] \ar[r] & G^G(X_1)
\oplus G^G(X_2) \ar[d] \ar[r]  & G^G(X_1 - X_1 \cap X_2) \oplus G^G(X_2
- X_1 \cap X_2) \ar[d]_{\simeq} \ar[r] & 0\\
G^G(X_1 \cap X_2) \ar[r] & G^G(X) \ar[r] &
G^G(X_1 - X_1 \cap X_2) \oplus G^G(X_2 - X_1 \cap X_2) \ar[r] & 0}
$$

A diagram-chase shows that the map $G^G(X_1) \oplus G^G(X_2)
\rightarrow G^G(X)$ is surjective.
\end{proof}

\begin{remark}
The analogous statement also
holds for non-equivariant higher $K$-theory. However, for higher
$K$-theory the only proof we know uses the Brown-Gersten
spectral sequence. Since we do not know how to adapt that
sequence to equivariant $K$-theory we can not state a
result for higher equivariant $K$-theory.
\end{remark}

\begin{lemma}\label{koszul} Let 
$L \subsetneq V$ be representations of $G$ and let 
$i: X \times L \rightarrow X \times V$ be the inclusion.  Then 
$$i_*G^G(X \times L) \subset I_G G^G(X \times V)$$
where $I_G \subset R(G)$ is the augmentation ideal.
\end{lemma}
\begin{proof}
Let $\pi_V: X \times V \rightarrow X$ and $\pi_L: X \times V
\rightarrow X$ be the projections.  By the homotopy property of
equivariant $G$-theory the smooth pullbacks $\pi_V^*$ and $\pi_L^*$
are isomorphisms of $R(G)$-modules. Since $i$ is a regular embedding
and $\pi_V \circ i = \pi_L$ we have, by the compatibility of flat and
l.c.i. pullbacks \cite[VI.6]{FL}, that $i^*\pi_V^* = \pi_L^*$.
Thus, the pullback $i^*:G^G(X \times V) \rightarrow G^G(X \times L)$
is an isomorphism of $R(G)$-modules.  Hence it suffices to show that
$i_*i^*: G^G(X \times L) \subset I_GG^G(X \times V)$. By the
projection formula
$i_*i^*$ acts by multiplication by $i_*(1) \in K^G(X \times V)$.
Since $L \subsetneq V$ is an invariant linear subspace, 
$i_*(1) = \lambda_{-1}(V/L)$ which is in $I_G$.
\end{proof}

\begin{lemma} \label{l.nilpotence}
Let $Z$ be a separated algebraic space, and let\\ $\a_Z = \ker (K(Z) 
\rightarrow \Z)$ be the augmentation  ideal. Then if $k > \mbox{dim}\;Z$,
$\a_Z^kG(Z) = 0$. If $Z$ is not separated then $\a_Z^kG(Z) = 0$
for $k$ sufficiently large.
\end{lemma}

\begin{proof}
Let $Z$ be a separated algebraic space. Following \cite[Definition
18.3]{Fulton} we define a Chow envelope $Z' \stackrel{p} \rightarrow
Z$ to be a proper morphism from a quasi-projective scheme $Z'$, such that
for every integral subspace $V \subset Z$, there is a subvariety $V'$
of $X$ such that $p$ maps $V'$ birationally to $V$.  Using Chow's lemma
for algebraic spaces \cite[Theorem IV.3.1]{knutson}, 
the argument of \cite[Lemma 18.3]{Fulton} shows
that every algebraic space has a Chow envelope (of the same
dimension as $Z$), and that the proper
push-forward $p_*:G(Z') \rightarrow G(Z)$ is surjective.

Since $Z'$ is quasi-projective, $F_{Z'}^k G(Z') = 0$
\cite[Cor. 3.10]{FL} where $F_{Z'}^k$ is the $k$-th
level of the $\gamma$-filtration. 
Since $p^*:K(Z) \rightarrow K(Z')$ preserves the rank
of a locally free sheaf, $p^*\a_Z^k \subset \a_{Z'}^k$. 
By definition, $\a_{Z'}^k \subset F^k$, so
$p^*\a_{Z}^k G(Z') = 0$. Hence, by the projection formula,
$\a_Z^k p_*G(Z') = 0$. Since $p_*$ is surjective, the lemma follows.

Even if $Z$ is not separated it has an open set $W$ which is a
separated scheme.  Then $\a^k G^G(W) = 0$ for $k >
\mbox{dim}\;Z$. Using the Noetherian induction and the localization
sequence we see that $\a_W^kG^G(Z) = 0$ for $k>>0$.
\end{proof}

\begin{lemma} \label{l.gnilp}
Suppose that $G$ acts freely on an algebraic space
$Z$. Then $I_G^kG^G(Z) = 0$ for $k>>0$.
\end{lemma}
\begin{proof}
Let $Y = Z/G$ be the quotient algebraic space. The proof
of \cite[Proposition 0.9]{GIT} extends to algebraic spaces and shows that
$Z \rightarrow Y$ is a $G$-principal bundle. Thus, there is an equivalence
between the categories of coherent sheaves on $Y$ and $G$-equivariant
sheaves on $X$; i.e., 
$G^G(Z) \simeq G(Y)$. Under this isomorphism, 
$I_G^lG^G(Z) \subset \a^lG(Y)$. 
The lemma now follows from Lemma \ref{l.nilpotence}.
\end{proof}

\begin{proof}[Proof of Theorem \ref{t.comp}]
We prove only part (a). The corresponding result about
Chow groups has essentially the same proof.

To show that the
filtrations of $G^G(X)$ by the submodules $\ker k_V$
and by powers of the ideal $I_G$ generate the same topology, there are
two steps.

{\it Step 1.}  We must show that given any pair $(V,U)$, there
is an integer $k$ 
such that $I_G^k G^G(X) \subset \ker k_V$, or in other
words, that $I_G^k G^G(X \times U) = 0$.  Since
$G$ acts freely on $X \times U$ this follows from Lemma \ref{l.gnilp}.

\medskip

{\it Step 2.} We must show that given a positive integer $k$, there is a
pair $(V,U)$ such that $\ker k_V \subset I_G^k G^G(X)$.
Let $(V,U)$ be any good pair and set $C = V-U$. 
Then $\ker k_V = \im(G^G(X \times C) \rightarrow G^G(X))$.
By definition, $C = \bigcup C_i$ where $C_i$ is contained
in an invariant linear subspace $L_i$ of $V$, and by Lemma \ref{koszul},
$\im(G^G(X \times C_i)) \subset I_GG^G(X)$.  By Lemma
\ref{component}, $\im(G^G(X \times C))$ is generated
by the images of $G^G(X \times C_i)$, so
$\im(G^G(X \times C)) \subset I_GG^G(X)$.

The group $G$ acts freely on the open set $U_k = V^{\oplus k} - C^k$, and
by induction,
$$\im(G^G(X \times C^k)) \subset I_G^kG^G(X \times V^{\oplus k}).$$
Thus, $(V^{\oplus k},U_k)$ is the desired good pair.
\end{proof}

\begin{example} \label{e.noncommute} There is a natural map
$$
\liminv (G^G(X)_{\Q}/\ker (k_V)_{\Q}) \rightarrow \liminv
(G^G(X)_{\Z}/\ker (k_V)_{\Z}) \otimes \Q
$$
However, because inverse limits do not commute with tensor product,
this map need not be an isomorphism.  
For example, if $G = {\mathbb G}_m$ is the 1-dimensional torus
and $X$ is a point, then 
$$
G^G(X)_{\Z} = K^G(X) = R(G) = \Z[u,u^{-1}].
$$
Let ${\mathcal V}$ denote the system $(V^{\oplus k}, U_k)$ where
$U_k = V^{\oplus k} - 0$.  Write \\
$x = u - 1$; then
$\liminv (G^G(X)_{\Z}/\ker (k_V)_{\Z})\simeq \Z[[x]]$,
so 
$$\liminv (G^G(X)_{\Z}/\ker (k_V)_{\Z})
\otimes \Q \simeq \Z[[x]] \otimes_{\Z}\Q.$$  This is
not isomorphic to
$$\liminv (G^G(X)_{\Q}/\ker (k_V)_{\Q}) \simeq \Q[[x]].$$
\end{example}

\section{Construction of an equivariant Riemann-Roch map}

If $Z$ is a separated algebraic
space, then there is a Riemann-Roch map 
$\tau_Z:G(Z) \rightarrow CH^*(Z)_{\Q}$ with the same properties
as the Riemann-Roch map for schemes constructed 
\cite[Theorem 18.3]{Fulton}. This fact \cite{gillet}
can be deduced from the Riemann-Roch theorem for
quasi-projective schemes and the existence of Chow envelopes
for separated algebraic spaces.

In this section we construct, for a separated $G$-space $X$, an equivariant Riemann-Roch map
$$\tau^G: G^G(X) \rightarrow \prod_{i=0}^\infty CH^i_G(X)_\Q$$
with the same functoriality as in the non-equivariant case 
\cite[Chapter 17]{Fulton}. In addition this map
will factor through the completion map $G^G(X) \rightarrow 
\widehat{G^G(X)}$.

The results of the previous section are not needed to construct
the map but they are used in the next section
to show that the map induces an isomorphism
$\widehat{G^G(X)}_{\Q} \rightarrow \prod_{i=0}^\infty CH^i_G(X)_\Q$.
To simplify notation we will write
$CH^i_G(X)$ for $CH^i_G(X)_\Q$.

All spaces in this section are assumed to be separated.

\subsection{Equivariant Todd classes and Chern characters.}
If $E \rightarrow X$ is an equivariant vector bundle of rank $r$
then it has equivariant Chern classes $c_1^G(E), \ldots , c_r^G(E)$
which are elements of the equivariant operational Chow ring $A^*_G(X)$
\cite[Sections 2.4, 2.6]{EIT}.
As in the non-equivariant case, an equivariant vector bundle $E \rightarrow
X$ of rank $r$ has Chern roots
$x_1, \ldots x_r$ such that
$c_i^G(E) = e_i(x_1, \ldots x_r)$
where $e_i(x_1, \ldots , x_r)$ refers to the $i$-th elementary
symmetric polynomial. 
This follows from the fact that the
equivariant Chern classes can be calculated on a fixed mixed space $X \times^G U$
\cite[Section 2.4, Definition 1]{EIT} and the non-equivariant
splitting principle \cite[Remark 3.2.3]{Fulton}.

\begin{defn} Let $E$ be a vector bundle with Chern roots
$x_1, \ldots x_r$.
Define the equivariant Chern character 
$$ch_G:K^G(X) \rightarrow \prod_{i = 0}^{\infty} A^i_G(X)$$
by the formula 
$$ch_G(E) = \sum_{i = 1}^{r} e^{x_i}.$$
Likewise, define the equivariant Todd class by the formula
$$\td^G(E) = \prod_{i= 1}^r \frac{x_i}{1-e^{-x_i}}.$$
\end{defn}
Because the leading coefficient of $\td^G(E)$ is $1$, 
$\td^G(E)$ is an invertible element of $\prod_{i = 0}^{\infty} A^i_G(X)$.  

\subsection{Construction of an equivariant Riemann-Roch map}
Recall that if $G$ acts freely on a space $Y$ then we have 
identifications
\begin{equation} \label{e.id}
G^G(Y) = G(Y/G) \mbox{  and  } CH^*_G(Y) = CH^*(Y/G).
\end{equation}
In what follows, we will use
such identifications, often without further comment, and when we
compare $\tau^G$ and $\tau_{Y/G}$ we are tacitly
using these identifications.

We want to define $\tau^G: G^G(X) \rightarrow \prod_{i=0}^\infty
CH^i_G(X)$ without assuming that the action is free, so that if the action 
is free then $\tau^G$
coincides with the non-equivariant map $\tau_{X/G}$.  We define
$\tau^G$ as follows.  
Choose a 
representation 
$V$ such that $G$ acts freely on an open set $U$ 
and $\mbox{codim}(V-U)>k$ (such pairs
always exist by \cite[Remark after Lemma 3]{EGchar}).  
The action of $G$ on $U$ and hence
on $X \times U$ is also free -- in particular, proper -- so
the quotient $X \times^G U$ is a separated algebraic space 
\cite[Corollary 2.2]{EdMo}.  Define
$$\rho_U: G^G(X \times U)\rightarrow CH^*_G(X \times U)$$ by
$$
\rho_U(\beta) = \frac{\tau_{X \times^G U}(\beta)}{\td^G(V)}.
$$

Under the identification $CH^*_G(X \times U) = CH^*(X \times^G U)$,
the action of $\td^G(V)$ is that of the Todd class of the vector
bundle $X \times^G(U \times V) \rightarrow X \times^G U$.

If $j < k$ we identify $CH^j_G(X)$ with $CH^j_G(X \times_G U)$.
If $\alpha \in G^G(X)$, we define $\tau^G(\alpha) \in
\prod_{i=0}^\infty CH^i_G(X)$ to be the element whose $j$-component
agrees with that of the image of $\alpha$ under the composition
$$
G^G(X) \rightarrow G^G(X \times U) \stackrel{\rho_U} \rightarrow CH^*_G(X \times U),
$$
where the first map is flat pullback.

\begin{prop} \label{p.taudef}  The definition of $\tau^G$ is
independent of the choice of $V$ and $U$.
\end{prop}

\begin{proof}
It suffices 
to show that
given $V \supset U$ and $V' \supset U'$ with $\mbox{codim}(V - U)$ and
$\mbox{codim}(V' - U')$ greater than $k$ the construction using
$V$ and $U$ agrees with that using $V'' = V \times V'$ and an open
subset $U''$ of $V''$.  We can choose the open subset of $V''$
arbitrarily, provided the codimension is sufficiently large, so we
take the open subset $U'' = U \times V'$.  The constructions agree
because the following diagram commutes.
$$
\begin{array}{ccc}
G(X \times^G(U \times V')) & \stackrel{\rho_{U \times V'} }
 \longrightarrow & CH^*(X \times^G(U \times V')) \\
\uparrow & & \uparrow\\
G(X \times^G U) & \stackrel{\rho_{U} }
 \longrightarrow & CH^*(X \times^G U)
\end{array}
$$
Here the vertical arrows are flat pullback, and commutativity
follows from \cite[Theorem 18.3(4)]{Fulton}, using the fact that the relative
tangent bundle of the morphism
\begin{equation} \label{e.morphism}
X \times^G(U \times V') \rightarrow X \times^G U
\end{equation}
is the bundle
$$
X \times^G(U \times V' \times V') \rightarrow X \times^G(U \times V').
$$
\end{proof}

\begin{remark} \label{r.vbs}
There is a variant construction of the map $\tau^G$,
which will be useful below.
Let 
${\mathcal V} \rightarrow X$ be an equivariant vector bundle 
such that $G$ acts freely
on an open set ${\mathcal U} \subset {\mathcal V}$ which surjects onto $X$.
If $\beta \in G^G(X \times {\mathcal U})$ set
$$\rho_{\mathcal U}(\beta) = 
\frac{\tau_{{\mathcal U}/G}(\beta)}{\td^G{\mathcal V}} \in CH^*_G(X \times 
{\mathcal U}).$$
If $\mbox{codim}({\mathcal V} - {\mathcal U})> k$, then the homotopy
property of equivariant Chow groups shows that we can identify
$CH^k_G({\mathcal U})$ with $CH^k_G(X)$ when $k < j$.
Arguing as in the proof of the proposition we can use pairs of
the form $({\mathcal V},{\mathcal U})$ to obtain
an equivalent definition of $\tau^G$.
\end{remark}

\begin{thm} \label{thm.tau}
Let $X$ be a separated $G$-space.
The equivariant Riemann-Roch map 
$$\tau^G: G^G(X) \rightarrow \prod_{i=0}^{\infty} CH^i_G(X)$$
has the following properties.

(a) $\tau^G$ factors through the completion map 
$G^G(X) \rightarrow \widehat{G^G(X)}$.

(b) $\tau^G$ is covariant for equivariant proper morphisms.

(c) If $\epsilon \in K^G(X)$ and $\alpha \in G^G(X)$
then $\tau^G(\epsilon \alpha) = ch^G(\epsilon) \tau^G(\alpha)$.

(d) Let $f:X \rightarrow Y$ be a morphism.
Assume either

(i) $f$ is a smooth and equivariantly quasi-projective\footnote{We 
say that a morphism $f:X
\rightarrow Y$ is equivariantly quasi-projective if there exists a
$G$-linearized line bundle on $X$ which is $f$-ample.  This is an
equivariant version of \cite[Definition 5.3.1]{EGA II}.}  morphism, or

(ii) $f:X \rightarrow Y$ is an
equivariant l.c.i. morphism,
$X$ and $Y$ can be equivariantly embedded in smooth $G$-schemes,
and $G$ is either special or connected.

Then $\tau^G(f^*\alpha) = \td^G(T_f) f^*\tau^G(\alpha)$,
where $T_f \in K^G(X)$ is the relative tangent element of the morphism $f$.

(e) If $G$ acts freely on $X$, then the map $\tau^G$ coincides with
the non-equivariant map $\tau_{X/G}$ under the identifications
$G^G(X) = G(X/G)$\\ and $CH^*_G(X) = CH^*(X/G)$.

\medskip

Moreover, $\tau$ is uniquely determined by properties (d(i)) and (e).
\end{thm}

\begin{proof}
Properties (b) and (c) follow from the definition of 
$\tau^G$ and
the non-equivariant Riemann-Roch theorem of \cite[Chapter 18]{Fulton}.
Property (e) follows because, as in Proposition \ref{p.taudef},
the diagram
$$
\begin{array}{ccc}
G(X \times^G V^0) & \stackrel{\rho_{V^0} }
 \longrightarrow & CH^*(X \times^G V^0) \\
\uparrow & & \uparrow\\
G(X/G) & \stackrel{\tau_{X/G} } \longrightarrow & CH^*(X/G)
\end{array}
$$
commutes.  

If $I_G$ is the augmentation ideal of $R(G)$, then
$ch^G(I_G^N) \subset \prod_N^\infty A^i_G$.
Thus, by property (c), 
$$\tau^G(I^N_G G^G(X)) \subset \prod_{N}^{\infty} CH^i_G(X).$$
Thus, $\tau^G$ restricts to a map 
$$G^G(X)/I_G^NG^G(X) \rightarrow \prod_{i=0}^N CH^i_G(X).$$
Taking the limit as $N \rightarrow \infty$ gives the desired
factorization
$$G^G(X) \rightarrow \widehat{G^G(X)} \rightarrow \prod_{i=0}^{\infty}
CH_G^i(X),$$
proving (a).

The proof of property (d(i)) uses \cite[Theorem 18.3(4)]{Fulton} but
requires some care, particularly in the category of algebraic spaces.
Suppose that $f:X \rightarrow Y$ is quasi-projective. Let $U$ be an
open set in a representation on which $G$ acts freely.  Denote the
mixed spaces $X \times^G U$ and $Y \times^G U$ by $X_G$ and $Y_G$.
Using descent \cite[Sections 8.4-8.5]{SGA} as in \cite[Proposition
2]{EIT}, we see that the induced map $X_G \rightarrow Y_G$ is smooth
and quasi-projective.  Let $Y' \stackrel{p} \rightarrow Y_G$ be a Chow
envelope for $Y_G$.  Then we have a Cartesian diagram
$$\begin{array}{ccc}
X' & \stackrel{f'} \rightarrow & Y'\\
q\downarrow & & p\downarrow\\
X_G & \stackrel{f} \rightarrow & Y_G\\
\end{array}$$
where, by base change,  $f'$ is smooth and quasi-projective.
Since $Y'$ is a quasi-projective scheme and $f'$ is a quasi-projective
morphism,
it follows from \cite[Corollary 5.3.3]{EGA II}
that $X'$ embeds in projective space as well. Hence
\cite[Theorem 18.3(4)]{Fulton} applies to the morphism
$f'$. If $\alpha \in G(Y_G)$ then $\alpha = p_*\alpha'$
for some $\alpha' \in G(Y')$.  Since flat pullback commutes with
proper pushforward, $f^*p_*\alpha' = q_*f^{'*}\alpha'$.  
Thus,

\begin{tabular}{ccll}
$\tau(f^*\alpha)$  & $=$ & $\tau(f^*p_*\alpha')$  &\\
& $=$ & $\tau(q_*f^{'*}\alpha')$ & \\
& $=$ & $q_*\tau(f^{'*}\alpha')$ & \cite[Theorem 18.3(2)]{Fulton}\\
& $=$ & $q_*[\td(T_{f'})f^{'*}\tau(\alpha')]$ 
& \cite[Theorem 18.3(4)]{Fulton}\\ 
& $=$ & $q_*[\td(q^*T_f)f^{'*}\tau(\alpha')]$ & (since $T_{f'} = q^*T_f$)\\
& $=$ & $\td(T_f)f^*\tau(\alpha)$ & (by the projection formula).
\end{tabular}

\noindent Taking the limit over all pairs $(V,U)$ in the construction of
$\tau^G$ gives (d(i)).

To prove (d(ii)) argue as follows: Suppose that $X \subset M$ and $Y
\subset Q$ where $M$ and $Q$ are smooth $G$-schemes. The requirement
that $G$ is special or connected ensures that we can choose open sets
$U \subset V$ so that $M \times^G U$ and $Q \times^G U$ are smooth
schemes \cite[Proposition 23]{EIT}. Thus the mixed space $X_G$ and
$Y_G$ are embeddable in smooth schemes, and again by descent
\cite[Proposition 2]{EIT}, the induced map $X_G \rightarrow Y_G$
is l.c.i., so (d(ii)) follows again from \cite[Theorem 18.3(4)]{Fulton}.

Finally, suppose that $\tau'$ is another map with
properties (d(i)) and (e). Suppose $\alpha \in G^G(X)$ is given
and denote by $\tau(\alpha)^j$ (resp. $\tau'(\alpha)^j$) the
term of $\tau(\alpha)$ (resp. $\tau'(\alpha)$) in $CH^j_G(X)$.
Let $(V,U)$ be a representation and open set on which $G$ acts
freely such that $\mbox{codim }(V-U) > j$. Let
$\pi_U: X \times U \rightarrow X\times V \rightarrow X$ be the composition
of projection with open inclusion.
Since $\mbox{codim }(V-U) > j$ we may 
identify $CH^j_G(X) = CH^j_G(X \times V) = CH^j_G(X \times U)$.
The morphism $\pi_U$ is smooth and quasi-projective so by property (d(i)) 
$$\tau'(\alpha)^j = \left(\frac{\tau'(\pi^*_U\alpha)}{\td^G(V)}\right)^j.$$
Since $G$ acts freely, $\tau'$ and $\tau$ coincide
on on $X \times U$.
Thus
$$\tau'(\alpha)^j = 
\left(\frac{\tau(\pi^*_U\alpha)}{\td^G(V)}\right)^j = \tau(\alpha)^j.$$
\end{proof}

\medskip

Let $E \rightarrow X$ be an equivariant vector bundle on 
a complete variety, and let $\pi$ denote the morphism $X \rightarrow pt$.  
Then $\pi_*(E) = \sum (-1)^i[H^i(X,E)] \in R(G)$.  
Set $\chi^G(E) = ch[\pi_*(E)] \in \prod_0^\infty A^i_G$, and
$\td^G(X) = \tau_X^G({\mathcal O}_X)$.
Applying the general
Riemann-Roch theorem to the morphism $X \stackrel{\pi} \rightarrow pt$
yields:
\begin{cor} \label{c.hirz}(Equivariant Hirzebruch-Riemann-Roch)
$$\chi^G(E)  = \pi_*(ch^G(E) \td^G(X))$$ in $\prod_0^\infty A^i_G$. \endproof
\end{cor}

\subsection{Example: The Weyl character formula} \label{s.weyl}
In this section we illustrate the equivariant Riemann-Roch theorem
by using it to prove the Weyl character formula for $SL_2$.
This is essentially a special case of a calculation done by
Bott \cite{Bott} using 
the Atiyah-Singer index theorem.  
This proof is different than some other proofs in that it does
not use a localization theorem to reduce to a computation
at the fixed point locus of the action of a maximal torus on the
flag variety. 

Let $T = \bG_m$ acting on $\P^1$ with weights $1,-1$.  What we will
calculate is $\chi^T({\mathcal O}_{\P^1}(n)) \in \prod_{i=0}^{\infty}
CH^i_T(pt)$.  This calculation is related to representations of
$SL_2$, since we can view $T$ as embedded in $SL_2$ as a maximal
torus.  Then $\sum (-1)^i H^i(\P^1, {\mathcal O}_{\P^1}(n))$ is a
virtual representation of $SL_2$ and $\chi^T({\mathcal O}_{\P^1}(n))$
gives a formula for the restriction of this representation to $T$.
We will prove that
$$\chi^T({\mathcal O}_{\P^1}(n)) = \frac{e^{(n+1)t} - e^{-(n+1)t}}
{e^t - e^{-t}}$$
(the notation is explained below).
For $n \geq 0$ the higher cohomology groups of ${\mathcal
O}_{\P^1}(n)$ vanish, the resulting representation of $SL_2$ is
irreducible, and this
gives the Weyl character formula.

Here is the calculation.  We have $A^*_T = \Q[t]$ \cite[Section 3]{EIT}, and 
$$\prod_{i=0}^{\infty}A^i_T
\cong \prod_{i=0}^{\infty} CH^i_T(pt)= \Q[[t]].$$
The Chern character $ch^T: R(T) \rightarrow \prod_{i=0}^{\infty}A^i_G(X)$
is given explicitly as follows.
If $V$ is a representation of $T$, write
$$
V = k_{n_1} \oplus \cdots \oplus k_{n_r},
$$
where $k_{n_i}$ is the $1$-dimensional representation of 
$T$ with weight $n_i$ (the weights $n_i$ need not be distinct).
Then $ch^T(V) = \sum e^{n_i t}$.

In \cite[Section 3.3]{EIT} we computed $A^*_T(\P^1) \simeq \Q[t,h]/(t +
h) (t -h)$ when $T$ acts with weights $\pm 1$. Under this
identification, $c_1^T({\mathcal O}_{\P^1}(n)) = nh$ and $c_1^T(T_{\P^1}) =
2h$.  Thus
$$ \pi_*(ch^T({\mathcal O}_{\P^1}(n)) \td^T(\P^1))= 
\pi_*(e^{nh} \frac{2h}{1- e^{-2h}}).
$$

\begin{lemma}
If $p(h) \in CH^*_T(\P^1) \simeq \Q[t,h]/(t + h) (t -h)$, then
$\pi_*(p(h))$ is given by
$$
\pi_*(p(h)) = \frac{p(t) - p(-t)}{2t}.
$$
The same formula holds for $\pi_*:\prod_{i=0}^{\infty} CH^i_T(\P^1) 
\rightarrow\prod_{i=0}^{\infty} CH^i_T(pt)$.
\end{lemma}

\begin{proof}
If $BT$ is a model for calculating $CH^*_T$, then the model
$\P^1_T$ for calculating $CH^*_T(\P^1)$ is a $\P^1$ bundle over
$BT$. Since $\pi_*(h) = 1$, the projection formula implies
that $\pi_*( \sum a_i t^i + h \sum b_j t^j)
= \sum b_j t^j$. Now if $p(h)  = \sum a_i h^i$ is in  $\prod_{i=0}^{\infty} CH^i_T(\P^1)$,
then we can write (using the relation $h^2 = t^2$) $x = \sum a_{2i} t^{2i}
+ h \sum a_{2i +1} t^{2i}$, so
$$
\pi_*(p(h)) = \sum a_{2i +1} t^{2i} = \frac{p(t) - p(-t)}{2t}
$$
which is the desired formula.
\end{proof}

The lemma implies that $\chi^T({\mathcal
O}_{\P^1}(n))$ is equal to
$$\pi_*(e^{nh} \frac{2h}{1 - e^{-2h}}) 
= \pi_*(2h \frac{e^{(n+1)h}}{e^h - e^{-h}}) = \frac{e^{(n+1)t} - e^{-(n+1)t}}
{e^t - e^{-t}}$$
as desired.

\subsection{Change of groups} \label{s.changegroups}
Let $H \subset G$ be an algebraic subgroup.
In
this section we discuss the relationship between $G$-equivariant
and $H$-equivariant Grothendieck groups, and the corresponding
Riemann-Roch maps $\tau^G$ and $\tau^H$.

\begin{lemma} \label{l.taucomp}
Given an action of $G \times H$ on $X$, with $H = 1 \times H$ acting
freely, there is a commutative diagram
$$
\begin{array}{ccc}
G^{G \times H}(X) & \stackrel{\tau^{G \times H}} \rightarrow 
& \prod CH^*_{G \times H}(X) \\
\simeq \downarrow & & \simeq \downarrow \\
G^G(X/H) & \stackrel{\tau^{G}} \rightarrow 
& \prod CH^*_G(X/H).
\end{array}
$$
\end{lemma}
\begin{proof}
Let $V$ be a representation of $G \times H$ and let $U$ be an
open set on which $G \times H$ acts freely. Then ${\mathcal V} = X \times^H V$
is a $G$-equivariant vector bundle on $X/H$. The group $G$
acts freely on the open set ${\mathcal U} = X \times^H U$ which surjects onto
$X/H$. Identifying  ${\mathcal U}/G$ with $X \times^{G \times H} U$
we see that the maps $\rho_U$ and $\rho_{{\mathcal U}}$
are the same. Using
pairs of the form $(V,U)$ we define $\tau^{G \times H}$ as $\rho_U$,
as before. 
On the other hand,
by Remark \ref{r.vbs} we can use
vector bundles and open sets of the form 
$({\mathcal V}, {\mathcal U}) = (X \times^H V,X \times^H U)$
to define $\tau^G_{X/H}$ as $\rho_{{\mathcal U}}$.  Hence
the maps coincide.
\end{proof}

\medskip

In $K$-theory, part (a) of the next proposition is due
to Thomason \cite{Tho87}, following ideas that go back to Atiyah and
Segal \cite{A-S}.  
Part (c) has apparently been used implicitly by Thomason \cite{Thoduke88}, 
but we do not know of an explicit statement or proof.

\begin{prop} \label{p.ind} Let $H$ be a subgroup of $G$, and let
$X$ be an $H$-space.

(a) There is a natural isomorphism of $R(G)$-modules $G^G(G \times^H
X) \simeq G^H(X)$, and a natural isomorphism of $A^*_G$-modules $\oplus
CH^i_G(G \times^H X) \simeq \oplus CH^i_H(X)$.

(b) The isomorphisms are compatible with the $\tau$ maps,
i.e., the following diagram commutes:
$$
\begin{array}{ccc}
G^{G}(G \times^H X) & \stackrel{\tau^{G}} \rightarrow 
& \prod CH^*_{G}(G \times^H X) \\
\cong \downarrow & & \cong \downarrow \\
G^H(X) & \stackrel{\tau^{H}} \rightarrow 
& \prod CH^*_H(X).
\end{array}
$$

(c) If the $H$-action on $X$ is the restriction of a $G$-action,
then $G \times^H X \cong G/H \times X$.  Under the isomorphism in (a),
the ``forgetful'' map $G^G(X) \rightarrow G^H(X)$ corresponds
to the flat pullback $G^G(X) \rightarrow G^G(G/H \times X)$.
A similar statement holds for Chow groups.
\end{prop}

\begin{proof}
We will only prove parts (a) and (c) for Grothendieck groups; the
arguments for Chow groups are similar.

Let $X$ be an $H$-space.
Define an $H \times G$-action
on $G \times X$ by:
\begin{equation}
\label{e.actionGtimesX}
(h,g) \cdot (g',x) = (g g' h^{-1}, h \cdot x)
\end{equation}
We also define an $H \times G$-action on $X$ by:
\begin{equation}
\label{e.action1}
(h,g) \cdot x = h \cdot x
\end{equation}
The projection $\pi: G \times X \rightarrow X$ is $H \times G$-equivariant, and
moreover a $G$-principal bundle ($G$-torsor). The projection 
$G \times X \rightarrow G \times^H X$ 
is $G$-equivariant and an $H$-torsor.  
Hence we have equivalences
of categories between:

$H$-equivariant coherent sheaves on $X$, 

$G \times H$-equivariant coherent sheaves on $G \times X$, and 

$G$-equivariant coherent sheaves on $G \times^H X$. \\ (This
uses the general fact \cite{Tho87} that if $M \rightarrow N$ is
a principal bundle, then there is an equivalence of categories
between coherent sheaves on $N$ and equivariant coherent sheaves on $M$.)  
Hence $G^G(G \times^H X) \simeq G^H(X)$.
Under these equivalences of categories, the 
$G$-equivariant vector bundle $V \times (G \times^H X) \rightarrow
G \times^H X$ corresponds to the $H$-equivariant vector bundle $V \times X
\rightarrow X$.  This translates into the fact that the isomorphism
is an isomorphism of $R(G)$-modules.  This proves (a) in $K$-theory. 
Part (b) follows from 
Lemma \ref{l.taucomp}.  

We prove (c).  We now assume that the $H$-action on $X$ is the restriction
of a $G$-action.  We can then define another $H \times G$-action on
$X$, this time by:
\begin{equation}\label{e.action2}
(h,g) \cdot x = g \cdot x.
\end{equation}
Let $\F$ be a $G$-equivariant sheaf on $X$.  Then $\F$ is naturally
an $H \times G$-equivariant sheaf, with respect to either action of $H \times G$.  
As noted above, the projection $\pi: G \times X \rightarrow X$ is
$H \times G$-equivariant with respect to action \eqref{e.action1} on $X$,
so $\pi^* \F$ is an $H \times G$-equivariant sheaf on $G \times X$.
Next, the action map $a: G \times X \rightarrow X$ is
$H \times G$-equivariant with respect to action \eqref{e.action2} on $X$.
Therefore, $a^* \F$ is also $H \times G$-equivariant.

In particular, equations \eqref{e.action1} and \eqref{e.action2} give
two $G \times G$-actions on $X$.  Since $\F$ is a $G$-equivariant
sheaf on $X$, by definition there is an isomorphism of coherent
sheaves $\theta: \pi^*\F \rightarrow a^*\F$ satisfying the cocycle
condition (see \cite[Section 1.2]{Tho87}).  The cocycle condition
implies:

\begin{lemma} \label{l.pi*a*} If $\F$ is a $G$-equivariant sheaf on
$X$, then the  map $\theta: \pi^* \F \rightarrow a^* \F$ is
an isomorphism of $G \times G$ equivariant sheaves on $G \times X$.  \endproof
\end{lemma}

Consider the diagram
$$
\begin{array}{ccccc}
G^G(X) & \stackrel{\pi^*}{\rightarrow} & G^{G \times G}(G \times X) & 
\rightarrow & G^G(G/G \times X) \\

\downarrow & & \downarrow & & \downarrow \\

G^H(X) & \stackrel{\pi^*}{\rightarrow} & G^{H \times G}(G \times X) & 
\rightarrow & G^G(G/H \times X ) 
\end{array}
$$
The left two vertical maps are the forgetful maps, and the
left square commutes.  The map $G^{H \times G}(G \times X)  
\rightarrow  G^G(G/H \times X )$ comes from taking the quotient
by $1 \times H$, and using the identification $(G \times X) / (1 \times H)$
with $G/H \times X$ taking $(g,x)$ to $(gH,gx)$.  Likewise, we
identify $G^{G \times G}(G \times X)$ with $G^G(G/G \times X) = G^G(X)$.  
Tracing through the definitions, the latter identification is given
by $(a^*)^{-1}$.  By Lemma \ref{l.pi*a*} the composition along the top row is
the identity. The right vertical arrow is flat pullback, and
the reader can verify that the right square also commutes.
This proves Proposition \ref{p.ind}.  \end{proof}

\section{The Riemann-Roch isomorphism}
All spaces considered here are again assumed to be separated, so that
we can apply the Riemann-Roch theorems of \cite{Fulton}.

Let $G \subset GL_n$ and let $I \subset
R(GL_n)$ be the augmentation ideal.
In this section we will use the notation
$\widehat{G^G(X)}$ to denote the $I$-adic
completion of $G^G(X) \otimes \Q$, and continue
to denote $CH^i_G(X) \otimes \Q$ by $CH^i_G(X)$.

The purpose of this section is to prove the following theorem.
\begin{thm} \label{thm.tauiso} Let $X$ be a separated algebraic space.
The map 
$$\tau^G:\widehat{G^G(X)} \rightarrow \prod_{i=0}^\infty CH^i_G(X)$$
is an isomorphism.
\end{thm}

\begin{remark}
By Corollary \ref{c.segal},  $\widehat{G^G(X)}$ is also
the completion of $G^G(X)$ with respect to $I_G \subset R(G)$.
\end{remark}

By Proposition \ref{p.ind} we obtain
a commutative diagram
$$\begin{array}{ccc}
\widehat{G^{GL_n}(GL_n \times^G X)} & \cong & \widehat{G^G(X)}\\
\downarrow \tau^{GL_n} &  & \downarrow \tau^G\\
\prod_{i=0}^{\infty}CH_{GL_n}^i(GL_n \times^G X) & \cong & 
\prod_{i=0}^{\infty}CH^i_G(X)
\end{array}
$$
Thus, to prove the theorem it suffices to show
$\tau^{GL_n}$ is an isomorphism.
We will do this in
two steps.

\subsection{$G = B$ is the group of upper triangular matrices}
Let ${\mathcal V} = \{V_n,U_n\}$ be the system of good representations and
open subsets constructed in the proof of Theorem \ref{t.good}.
Since $B$ acts freely on $X \times U_n$, Theorem \ref{thm.tau}(e)
and the non-equivariant Riemann-Roch theorem
(\cite[Theorem 18.3]{Fulton}, extended to separated algebraic spaces in \cite{gillet})
implies that
$$\tau^B:G^B(X \times U_n) \rightarrow 
CH^*_B(X \times U_n)$$
is an isomorphism.

By Theorem \ref{thm.tau}, the map $\tau^G$ is compatible with the
maps in the inverse system of representations
constructed in Section \ref{s.completions}. 
In this way we obtain an isomorphism
$$\tau^B: \liminv (G^B(X)/\ker k_n) \rightarrow 
\liminv CH^*_B(X)/\ker r_n.$$

By Proposition \ref{p.borelsegal}, the $I_B$-adic completion
of $G^B(X)$ is isomorphic to the $I$-adic completion $\widehat{G^B(X)}$,
so by Theorem \ref{t.comp},
$$\liminv (G^B(X)/\ker k_n) \simeq \widehat{G^B(X)}.$$
Also, by Theorem \ref{t.comp} and Proposition \ref{p.brion},
$$\liminv CH^*_B(X)/\ker r_n \simeq \widehat{CH^*_B(X)}
 \simeq \prod_{i=0}^{\infty} CH^i_B(X).$$
The case $G=B$ follows.

\subsection{$G = GL_n$}
Let $\pi: G/B \times X \rightarrow X$ be the projection to the second
factor.  By Proposition \ref{p.ind}, we identify $G^B(X)$ with $G^G(G/B \times
X)$ so that the forgetful map $i^!:G^G(X) \rightarrow G^B(X)$
corresponds to flat pullback $\pi^*:G^G(X) \rightarrow G^G(G/B \times
X)$.

Since $B$ is the group of upper triangular
matrices, the quotient $G/B$ is complete and
we can define a map $i_!:G^B(X) \rightarrow G^G(X)$
by using the identification $G^B(X) \cong
G^G(G/B \times X)$ and setting $i_!$ to be the proper push-forward $\pi_!:G^G(G/B \times X)
\rightarrow G^G(X)$.

\begin{lemma} The composition $i_!i^! = id$, so
we can identify $G^G(X)$ as a summand in $G^B(X)$.
\end{lemma}
\begin{proof}
This follows from the argument of \cite{A-S}, using the fact that
$G/B$ is a tower of projective bundles and the projective bundle
theorem of \cite[Theorem 3.1]{Tho87}.
\end{proof}

\medskip

We now turn to the Chow theory.  
We use the identification
$CH^*_B(X) \simeq CH^*_G(G/B \times X)$
to define a map
$$i^*:\prod_{0}^{\infty}CH_G^i(X)  \rightarrow \prod_{i=0}^{\infty}
CH_B^i(X) \simeq
\prod_{i=0}^{\infty} CH_G^i(G/B \times X)$$ by the formula
$$
\alpha \mapsto \td^G(T_{G/B}) \pi^*(\alpha).
$$
We also set
$$
i_*: CH^*_B(X) \rightarrow CH^*_G(X)
$$
to be the composition
$$CH^*_B(X) \simeq CH_G^*(G/B \times X) 
\stackrel{\pi_*} \rightarrow CH_G^*(X)
.$$
\begin{lemma}
The composition $i_*i^*=id$, so we can identify $CH^*_G(X)$ as a summand 
in $CH^*_B(X)$.
\end{lemma}
\begin{proof}
In general, if $f: M \rightarrow N$ is a smooth proper morphism
with $f_! [\O_M] = [\O_N]$, then the Riemann-Roch theorem implies that\\
$f_*(\td(T_f) f^* \beta) = \beta$ for all $\beta \in CH^*(N)$.  
(Proof: By the projection formula we may assume $\beta =1$.
Then by Riemann-Roch $f_*(\td(T_f)) = 1$.)
Applying
this to the case where $M$ and $N$ are the mixed spaces
$(G/B \times X)\times^G U$ and $X \times^G U$ yields the result.
\end{proof}

\begin{lemma} \label{l.gmodb}
We have:

(i) $i_!\tau^G = \tau^Bi_*$

(ii) $\tau^Bi^! = i^*\tau^G$.
\end{lemma}
\begin{proof}
In view of the previous lemmas (i) follows from Theorem \ref{thm.tau}(a) and
(ii) from 
Theorem \ref{thm.tau}(d(i))
\end{proof}

\medskip

Since $\tau_B$ is an isomorphism the previous lemmas imply
$\tau_G$ is also an isomorphism when $G = GL_n$.  This proves
the theorem.
\endproof

\section{Actions with finite stabilizers} \label{s.finite}
In this section we assume that the group $G$ acts on the space
$X$ with finite stabilizers and rational coefficients are
used throughout.

In this situation, we can obtain
more refined information about the completions in Section \ref{s.completions}.
In particular we prove the following theorem.

\begin{thm} \label{t.locs}
Suppose that an algebraic group $G$ acts on a space
$X$ with finite stabilizers. Let $I_G$ be the augmentation ideal.
Then there is an isomorphism
$$\widehat{G^G(X)}  \rightarrow G^G(X)_{I_G}.$$
\end{thm}

We begin with a commutative algebra lemma.

\begin{lemma} \label{l.module1}
Let $R$ be a ring and $I$ a maximal ideal and let $M$ be an
$R$-module.  Suppose that $I^k M_I = 0$.  
Then $M_I = \widehat{M}$
where $\widehat{\mbox{  }}$ denotes the $I$-adic completion.
\end{lemma}
\begin{proof}
The rings $R/I^k$ are local, so $(R/I^k)_I = R/I^k$.
Thus, 
$$M_I/I^kM_I := M \otimes (R_I/I^k R_I) = M \otimes R/I^k := M/I^k M.$$
Hence $\widehat{M_I} \cong \widehat{M}$.  Since $I^k M_I = 0$,
$M_I$ is $I$-adically complete, proving the result.
\end{proof}

\medskip

Because of the lemma, Theorem \ref{t.locs} is an immediate
consequence of the following result, whose proof uses
Corollary \ref{c.segal}.

\begin{prop} \label{p.gsupp}
If $G$ acts on $X$ with finite stabilizers then for sufficiently large $k$,
$$I_G^kG^G(X)_{I_G} = 0.$$
\end{prop}

\begin{proof} 
The proof proceeds as in previous proofs, by building up from 
a torus to $GL_n$ to a general group.

\medskip

{\it Step 1.  $G = T$ is a torus.}
By Thomason's generic slice theorem \cite[Proposition 4.10]{Thoduke86}, 
there is a $T$-equivariant
open subset $U \subset X$ such that $U$ is equivariantly
isomorphic to $S \times Y$ where $S = T/H$ for a
diagonalizable subgroup $H \subset T$, and $T$ acts trivially on 
$Y$.  If $W$ is any representation of $T$ which is trivial on $H$,
then the vector bundle  $U \times W \rightarrow U$ is equivariantly
trivial, so $([W] - \mbox{dim }W) G^T(U) = 0$.  Hence, under
the map $R(S) \rightarrow R(T)$, the ideal $J = I_S R(T)$ annihilates
$G^T(U)$.  Note that $J \subset I_T$.

Since the stabilizers are finite $H$ is finite.
Let $\widehat{T}$ be the character group of $T$,
and let $\widehat{S}$ be the character group of $S$. Then we 
can choose a basis $e_1,\ldots,e_n$ of
$\widehat{T}$ such that $d_1 e_1,\ldots,d_n e_n$ is a basis
for $\widehat{S}$, where the $d_i$ are positive integers.  
Using this basis, we have 
$$R(T) \cong \Q[\widehat{T}]
\cong \Q[t_1,\ldots,t_n]$$ and 
$$J = (t_1^{d_1} -1, \ldots, t_n^{d_n} - 1) 
\subset I_T = (t_1 -1, \ldots,t_n - 1).$$  
Hence $J_{I_T} = I_T R(T)_{I_T}$, so
$I_T G^T(U)_{I_T} = 0$. 

By Noetherian induction we may assume that $I_T^kG^T(X-U)_{I_T} = 0$
for some sufficiently large $k$.  Then by the localization
exact sequence, $I_T^{k+1}G^T(X)_{I_T} = 0$.

Note that the above discussion and Noetherian induction also show that there
is a non-zero ideal $J_T \subset I_T$ such $J_TG^T(X) = 0$.
\medskip

{\it Step 2. $G = GL_n$.}  Let $B$ be the group of upper triangular
matrices in $GL_n$ and let $T$ be the group of diagonal matrices.
Restriction from $B$ to $T$ induces isomorphisms $R(B) \cong R(T)$ and
$G^B(X) \cong G^T(X)$ (see \cite[Proof of Theorem 1.11]{Thoduke88}).  
Claim: There is an ideal $J_G \subset I_G$
with the support of $R(G)/J_G$ finite, such that $J_G G^G(X) = 0$.  For,
$R(G) \cong R(B)^W \hookrightarrow R(B)$ is finite, so if
$J_G$ is the inverse image of $J_B := J_T$ under this map, then
the support of $R(G)/J_G$ is finite.  Since $G^G(X)$ is an $R(G)$-submodule
of $G^B(X)$ (see the proof of Theorem \ref{thm.tauiso}), we have
$J_G G^G(X) = 0$, proving the claim.

Because $R(G)/J_G$ has finite support and $J_G$ is contained in the
maximal ideal $I_G$, it follows that $J_G R(G)_{I_G}$ is
$I_G$-primary, hence contains $I_G^k R(G)_{I_G}$ for some $k$ (note
that $R(GL_n) = \Q[t_1,\ldots , t_n]^{S_n}$ is  Noetherian).  
Then $I_G^kG^G(X)_{I_G} = 0$ as desired.

{\it Step 3. The general case.}
Embed $G \subset GL_n$.  In this case, 
$G^G(X) = G^{GL_n}(GL_n \times^G X)$. Thus, by Step 2, 
$I_{GL_n}^k G^G(X)_{I_{GL_n}} = 0$ for some positive integer $k$. Equivalently, given
$x \in G^G(X)$ and $a \in I_{GL_n}^k$, there exists $b \in R(GL_n) - I_{GL_n}$
such that $abx = 0$.

Now, the action of $R(GL_n)$ on $G^G(X)$ factors through $\phi:R(GL_n)
\rightarrow R(G)$.  Also, the augmentation map $R(GL_n) \rightarrow \Q$
factors through $\phi$, so if $x \in R(GL_n) - I_{GL_n}$ then $\phi(x) \in
R(G) - I_G$.  Corollary \ref{c.kock} implies that for some
$d$, $I_G^d \subset \phi(I_{GL_n}) R(G)$.  This implies that
$I_G^{dk} G^G(X)_{I_G} = 0$.  Indeed, if $a' \in I_G^{dk}$ and
$x \in G^G(X)$, write $a' = \phi(a)$ for $a \in I_{GL_n}^k$,
choose $b$ as in the preceding paragraph; then $b' = \phi(b)
\in R(G) - I(G)$ and $a' b' x = 0$.  This proves the result.
\end{proof}

\medskip

\begin{remark} \label{r.annihilate}
The proof of the proposition implies that for $G=T$ or $G=GL_n$, that
there is an ideal $J \subset I_G \subset R(G)$ with the support of
$R(G)/J$ finite, such that $JG^G(X) = 0$. 

If $G$ is a diagonalizable group, embedded
in a torus $T$, then we know from the general theory of characters
of diagonalizable groups \cite[Section 8.12]{Borel},
that the map $R(T) \rightarrow R(G)$ is finite. Thus, an ideal
$J$ exists as above that annihilates $G^G(X)$.

In characteristic 0, the
 the map $R(GL_n) \rightarrow R(G)$ is 
finite for arbitrary $G$ \cite[Proposition 3.2]{Segal}, and
we can make the same conclusion about the support of $G^G(X)$. 
However, in characteristic $p$ we don't even know
if $R(G)$ is always Noetherian!
\end{remark}

\medskip

\begin{prop}
If $G$ acts with finite stabilizers then the rational equivariant
Chow groups $CH^*_G(X)_{\Q}$ are generated by invariant cycles
on $X$.  In particular $CH^i_G(X)_{\Q} = 0$ for 
$i > \mbox{dim }X$.
\end{prop}

\begin{proof}
By the generalization of \cite[Theorem 6.1]{Se} to algebraic spaces
there is a finite
cover $f:X' \rightarrow X$ on which $G$ acts freely.
Since $G$ acts freely, $CH^*_G(X')$ is generated 
by invariant cycles \cite[Proposition 8]{EIT}.
On the other hand, the proper pushforward
$f_*: CH^*_G(X')_\Q \rightarrow CH^*_G(X)_\Q$ is surjective
because $f$ is finite and surjective. Therefore, 
$CH^*_G(X)_\Q$ is generated by invariant cycles.
\end{proof}

\medskip

The preceding proposition implies that the equivariant Chow groups are complete.
Using the fact that $I$ and $I_G$ generate the same topology (Corollary 
\ref{c.segal}) we can restate the Riemann-Roch theorem.

\begin{cor} \label{c.fs}
Suppose that $X$ is separated and $G$ acts with finite stabilizers.
There is a map
$$\tau:G^G(X) \rightarrow G^G(X)_{I_G} \stackrel{\sim} \rightarrow 
CH^*_G(X)_\Q$$
satisfying properties (a)-(e) of Theorem \ref{thm.tau}.  \endproof
\end{cor}

\begin{remark}
When $G$ acts on properly on a separated scheme $X$ 
with {\it reduced} stabilizers then Vistoli
\cite{Vi3} stated a theorem which asserted the existence of
a map 
$$\t_X: G^G(X) \otimes \Q \rightarrow CH^*([X/G] \otimes \Q)$$
where here $[X/G]$ is the
Deligne-Mumford quotient stack.  By \cite[Proposition 14]{EIT}, 
$$CH^*([X/G])
\otimes \Q= CH^*_G(X) \otimes \Q.$$
Vistoli noted that his 
map need not be an isomorphism and made a conjecture 
about its kernel \cite[Conjecture 2.4]{Vi3}. The conjecture
states that if $\alpha \in \ker \tau_X$ then $\xi \alpha = 0$
for some $\xi$ which is the class of a perfect complex with everywhere
non-zero rank.

We expect that his map is the same as ours in this case.
Unfortunately, because he did not write his map down,
and he did not state whether it satisfied properties
d(i) and (e) of Theorem \ref{thm.tau},
we can not positively assert this.

However, for our map $\tau_X$ Vistoli's conjecture is true and 
his statement can be refined.
\begin{cor}(Vistoli's conjecture) \label{c.vistoli}
$$\alpha \in \ker(\tau^G_X:G^G(X) \rightarrow CH^*([X/G])_{\Q})$$
if and only there exists a virtual representation $\epsilon \in R(G)$ of
non-zero rank such that
$\epsilon \alpha = 0$.
\end{cor}
\begin{proof} The kernel is just the kernel of the localization
map $G^G(X) \rightarrow G^G(X)_{I_G}$.
\end{proof}
\end{remark}

\subsection{The case $G$ is diagonalizable} 

By Remark \ref{r.annihilate}, $G^G(X)$ is supported as an $R(G)$-module
at a finite number of primes, each of which is maximal.
Denote these ideals by $P=P_0, P_1, \ldots P_k$.
Following \cite{Segal},
each prime $P_i$ corresponds to a finite
subgroup (called the support of $P_i$) $H_i \subset G$. It is defined
as the minimal element of the set of subgroups $H \subset G$
such that $P_i \in Im(\mbox{Spec }R(H) \rightarrow \mbox{Spec }R(G))$.
Note that different $P_i$'s may have the same support.
This definition makes sense for any group $G$, but $H_i$ is only defined
up to conjugation as a subgroup of $G$.  In our case $G$ is abelian so the $H_i$'s
are uniquely determined. 

Following \cite[Lemma 1.1, Proposition 1.2]{Thoduke92}
we give an explicit construction of the support $H$ of
a prime ideal $P \subset R(G)$.
Since $G$ is diagonalizable,  $R(G) = \Z[N]$ where $N$ is a finitely
generated abelian group without $p$-torsion (where $p = \mbox{char }k$).
Given a prime $P \subset R(G)$ 
set
$$K_P = \{n \in N| 1 -n \in P\}.$$
The equivalence of categories between finitely generated abelian
groups (without $p$-torsion) and diagonalizable groups
means that quotient $N/K_P$ determines
a unique subgroup $H \subset G$ with the property that $R(H) = \Z[N/K_P]$.
When $P$ is maximal, $K_P$ has finite index in $N$ and
$H$ is a finite group.

The representation ring of the quotient $G/H_i$ is the subring
$\Z[K_P]$ of $\Z[N] = R(G)$. If $\Ibar$ is the augmentation ideal
of $R(G/H)$ then this construction shows that $\Ibar = P \cap R(G/H)$.

Denote by $X^i$ the subscheme fixed by $H_i$.
\begin{thm} \label{t.akt} Suppose that $G$ is diagonalizable
and acts on a separated space $X$ with finite stabilizers.
Then
$$G^G(X)_{P_i} \simeq CH^*_{G/H_i}(X^i) \otimes_{R(G/H_i)} R(G)_{P_i}.$$
In particular, there is an isomorphism 
$$G^G(X) \simeq \prod_{i} CH^*_{G/H_i}(X^i) 
\otimes_{R(G/H_i)} R(G)_{P_i}.$$
\end{thm}

\medskip


\medskip

\begin{proof}[Proof of Theorem \ref{t.akt}]

By the localization theorem for diagonalizable group
schemes $G^G(X)_{P_i} \simeq G^G(X^i)_{P_i}$ \cite[Theorem 2.1]{Thoduke92}.
Since $H_i \subset G$ acts trivially on $X^i$ we have \cite[
Lemma 5.6]{Thoinv}
$$ G^G(X^i) \simeq G^{G/H_i}(X^i) \otimes_{R(G/H_i)} R(G).$$
Let $\Ibar_i$ be the augmentation ideal of $G/H_i$. 
Since $P_i \cap R(G/H_i) = \Ibar_i$,
$$[G^{G/H_i}(X^i) \otimes_{R(G/H_i)} R(G)]_{P_i}
= [G^{G/H_i}(X^i)_{\Ibar_i} \otimes_{R(G/H_i)} R(G)]_{P_i}.$$
Thus, by the Riemann-Roch isomorphism of Corollary \ref{c.fs}
we have 
$$[G^{G/H_i}(X^i)_{\Ibar_i} \otimes_{R(G/H_i)} R_G]_{P_i}
\simeq [CH^*_{G/H_i}(X^i) \otimes_{R(G/H_i)} R_G]_{P_i}.$$
(Here the Chern character $ch: R(G/H_i) \rightarrow
\prod_{i=0}^{\infty} A^*_{G/H_i}$ makes $CH^*_{G/H_i}(X^i)$ into an $R(G/H_i)$-module.)
The first statement follows.

As noted in Remark \ref{r.annihilate} there is an ideal
$J \subset R(G)$ such that $R(G)/J$ is supported at a finite number of points
and
$JG^G(X) = 0$. Hence, $G^G(X) \simeq G^G(X) \otimes_{R(G)} R(G)/J$.
Then $J = Q_1 \cap  Q_2 \ldots \cap Q_k$ with $Q_i$ a
$P_i$-primary ideal. By the Chinese remainder theorem 
$$R(G)/J \simeq \prod_{i=1}^{r} R(G)/Q_i = \prod_{i=1}^{r} (R(G)/J)_{P_i}$$
so 
$$G^G(X) \simeq \prod_{i} G^G(X)_{P_i}$$
and the second statement follows.
\end{proof}

\begin{remark} This result was first obtained
by Angelo Vistoli (unpublished). It is also
related to the Riemann-Roch theorem for algebraic stacks
proved by B. Toen \cite{Toen}.
\end{remark}

\section{More on completions}
There are other natural completions of $G^G(X)$ and
$CH^*_G(X)$ besides those of Section \ref{s.completions}.
The purpose of this section is to prove that the different
definitions give isomorphic completions.  As an application
of these results, we will prove a special case of a conjecture of
K\"ock.  To begin, fix an embedding of $G$ into $GL_n$.
Then $G^G(X)$ is an $R(GL_n)$-module.

Recall that $\widehat{G^G(X)}$ was defined to be completion
of $G^G(X)$ along the augmentation ideal $I$ of $R(GL_n)$, and
$\widehat{CH^*_G(X)}$ the completion of $CH^*_G(X)$ along the
augmentation ideal of $A^*_G(pt)$.  We will refer to these as
the ``point'' completions because they are defined using ideals
in the equivariant groups of a point.

Let $I_X \subset K^G(X)$ denote the augmentation ideal, i.e., the
ideal of virtual vector bundles of rank $0$, and let
$\widetilde{G^G(X)}$ denote the completion of $G^G(X)$ along $I_X$.
Let $J_X \subset A^*_G(X)$ denote the augmentation ideal, and
$\widetilde{CH^*_G(X)}$ the completion of $CH^*_G(X)$ along $J_X$.  We
will refer to these as the ``$X$'' completions because they are
defined using ideals in the equivariant groups of $X$.  In the proof
below it will be necessary to distinguish between ideals corresponding
to different groups.  We will use a subscript to indicate this, e.g.,
$I_{X,G} \subset K^G(X)$.

Note that the $X$-completions only depend on $G$ and $X$, while
a priori, the point completions depend in addition on
an embedding of $G$ into $GL_n$.  Corollary \ref{c.segal} will imply
that the point completion is independent of the embedding.

The main result of this section is that the point and $X$-completions
are isomorphic.

\begin{thm} \label{t.compiso}

(a) The $I$-adic and $I_X$-adic topologies on $G^G(X)$ coincide.  
Hence we have an isomorphism of completions
$$\widetilde{G^G(X)} \simeq
\widehat{G^G(X)}.$$

(b) The $J$-adic and $J_X$-adic topologies on $CH^*_G(X)$ coincide.  
Hence we have an isomorphism of completions
$$\widetilde{CH^*_G(X)} \simeq \widehat{CH^*_G(X)}.$$
\end{thm}

\begin{proof}
(a) To show that the filtrations induced by powers of the
ideals $I$ and $I_X$ induce the same topology, we must
check two things.  First, we must show that for any $n$, there
exists an $r$ such that $I^r G^G(X) \subseteq I^n_X G^G(X)$.  
This is clear because under the map
$R(GL_n) \rightarrow K^G(X)$, the image of
$I$ is contained in $I_X$, so we can take $n = r$.
Second, we must show that for any $n$, there
exists an $r$ such that $I^r_X G^G(X) \subseteq  I^n G^G(X)$.  
As above we will do this in steps: first for $G=B$ the group 
of upper triangular matrices,
then for $G=GL_n$, and finally for arbitrary $G$.

Suppose then that $B \subset GL_n$ is the group of upper
triangular matrices.  We use the notation of the proof of
Theorem \ref{t.good}.  We know there exists $m$ such
that $\ker k_m \subseteq  I_B^n G^B(X)$.  Since
$I_B$ and $I$ generate the same topology on $R(B)$,
we can assume
$\ker k_m \subseteq I^n G^B(X)$. So we must
show that there exists $r$ such that 
$I^r_X G^B(X) \subseteq  \ker k_m$,
i.e., such that $I_X^r G^B(X \times U_m) = 0$.  Since $B$ acts
freely on $X \times U_m$, we have $G^B(X \times U_m) \simeq
G(X \times^B U_m)$.  Under this isomorphism, 
$$I_X^r G^B(X \times U_m) \subseteq \a^r G(X \times^B U_m)$$ 
where $\a^r$ denotes is
the augmentation ideal of $K(X \times^B U_m)$.
By Lemma \ref{l.nilpotence}, $\a^r G(X \times^B U_m) = 0$
for $r >>0$.

The analogous statement for Chow groups, that
$J^r_X CH^*_B(X \times U_m) = 0$ for $r > \mbox{dim}(X \times^B U_m)$,
also holds.

\medskip

Assume now that $G = GL_n$.  Then we have
$$
G^B(X) \simeq G^G(G/B \times X) \stackrel{\leftarrow} \rightarrow G^G(X).
$$
where the two maps are $i_!$ and $i^!$.
We have proved that there exists $r$ such that
$$
I^r_{X,B} G^B(X) \subset I^n G^B(X).
$$
Hence
$$
i_! (I^r_{X,B} G^B(X) ) \subset i_! (I^n G^B(X) ) =  I^n G^G(X)
$$
where the last equality follows by the projection formula.   
Now, we also have $i^! (I^r_{X,G} G^G(X) ) \subset
I^r_{X,B} G^B(X)$.  Combining these facts, we see that
$$
i_! i^! ( I^r_{X,G} G^G(X) ) \subset I^n G^G(X).
$$
Since $i_! i^!$ is the identity
$I^r_{X,G} G^G(X)  \subset I^n G^G(X)$.

\medskip

Finally consider the case where $G \hookrightarrow GL_n$ is any
subgroup.  By Proposition \ref {p.ind}, there is an isomorphism
or $R(G)$-modules
$$
G^G(X) \simeq G^{GL_n}(GL_n \times^G X).
$$
Under this isomorphism 
$I^n G^G(X)$ corresponds to $I^n G^{GL_n}(GL_n \times^G X)$.  
Moreover, 
$I^r_{X,G}  G^G(X)$ corresponds to 
$I^r_{GL_n \times^G X, GL_n} G^{GL_n}(GL_n \times^G X)$ since
the equivalence
of categories, obtained from descent of $GL_n \times G$-equivariant
sheaves on $GL_n \times X$, 
takes locally free sheaves to locally free
sheaves of the same rank \cite[Prop. 2.5.2]{EGA IV}.]  Because
of these correspondences,
the theorem follows from the case where $G = GL_n$.
\end{proof}

\medskip

\begin{cor} \label{c.segal}
If $H \subset G$, then the topology on $R(H)$ induced
by the ideal $I_{G}R(H)$ is the same as the $I_H$-adic topology.
\end{cor}
\begin{proof}
By embedding $H \subset G \subset GL_n$ it suffices
to prove the result for $G \subset GL_n$. The corollary
now follows by applying the theorem when $X = pt$.
\end{proof}

\begin{remark} In characteristic 0 this is the same as
\cite[Corollary 3.9]{Segal}, but in characteristic $p$ the
result is new.
For a large class of group schemes over perfect fields
Thomason \cite[Corollary 3.3]{Thoduke86} 
showed that $I_G$-adic and $I_{GL_n}$-adic topologies
are the same on the $\mbox{mod } {l}^{\nu}$
equivariant $G$-theory localized at the Bott element.
\end{remark}

Observe that the higher equivariant $K$-groups $G_i^G(X)$ (resp.
$K_i^G(X)$) are also modules over $R(G)$ and $K^G(X)$.
Thus we can define completions of these groups with respect
to the ideals $I_X$ and $I_G$.
If $X$ is a regular scheme then $K_i^G(X) = G_i^G(X)$ 
and we can prove a corollary about higher
$K$-theory as well. 
Part (b) proves a conjecture
of K\"ock \cite[Conjecture 5.6]{kock} for regular schemes over fields.
\begin{cor} \label{c.kock}
(a) If $X$ is a regular $G$-scheme then 
$\widetilde{K_i^G(X)}  \simeq \widehat{K_i^G(X)}$.

(b) Let $f:Y \rightarrow X$ be an equivariant proper morphism of regular 
$G$-schemes.
Then the push-forward $f_*: K_i^G(Y) \rightarrow K_i^G(X)$ induces a map
of completions
$f_*: \widetilde{K_i^G(Y)} \rightarrow \widetilde{K_i^G(X)}$.
\end{cor}

\begin{proof}
The action of $R(G)$ on $K_i^G(X)$ factors through the
map $R(G) \rightarrow K^G(X)$. 
Since $K^G(X) = G^G(X)$, Theorem \ref{t.compiso} implies that
the ideals $I_X \subset K^G(X)$ and the $I_GK^G(X)$
generate the same topology.
Thus, 
$\widetilde{K_i^G(X)}  \simeq \widehat{K_i^G(X)}$ proving (a).

By the projection formula applied to the commutative triangle
$$\begin{array}{ccc}
Y & \rightarrow & X\\
\downarrow &  \swarrow & \\
pt & & 
\end{array}
$$
$f_*(I^kK_i^G(Y)) = I^k f_*K_i^G(Y) \subset I^kK_i^G(X)$.
Hence, $f_*$ is continuous with respect to the $I$-adic topology, proving
(b).
\end{proof}

\begin{remark}
In its full form, K\"ock's conjecture asserts that if $G/S$ is
a flat group scheme and
if $X \rightarrow Y$ is any equivariant projective local complete
intersection morphism then there is a push-forward 
$$f_*:\widetilde{K_i^G(X)} \rightarrow \widetilde{K_i^G(Y)}$$
of completions. 
This conjecture is quite subtle because (despite the suggestive
notation) the completions
are taken with respect to different ideals, and if $X$ and $Y$
are not regular, there is no obvious way of comparing the topologies.
\end{remark}
\begin{remark} The $K_0$ version of K\"ock's conjecture has been proved
\cite{CEPT}
for  finite group schemes acting on regular projective varieties over
rings of integers of number fields.
\end{remark}

\end{document}